
\input eplain.tex

\vsize=21cm 
\voffset=1.5cm
\baselineskip=13pt 
\footline={\hss{\vbox to 1.5cm{\vfil\hbox{\rm\folio}}}\hss}

\input amssym.def
\input amssym.tex

\font\small=cmr8
\font\csc=cmcsc10
\font\title=cmbx10
\font\stitle=cmmib10

\font\csc=cmcsc10
\font\title=cmbx10 at 12pt
\font\stitle=cmmib10 at 12pt

\font\teneusm=eusm10
\font\seveneusm=eusm7
\font\fiveeusm=eusm5
\newfam\eusmfam
\def\eusm{\fam\eusmfam\teneusm}
\textfont\eusmfam=\teneusm
\scriptfont\eusmfam=\seveneusm
\scriptscriptfont\eusmfam=\fiveeusm

\def\varGamma{{\mit \Gamma}}
\def\Re{{\rm Re}\,}
\def\Im{{\rm Im}\,}

\def\scr#1{{\scriptstyle{#1}}}
\def\r#1{{\rm #1}}
\def\B#1{{\Bbb #1}}
\def\e#1{{\eusm #1}}

\def\sgn{{\rm sgn}}

\singlecolumn
\centerline{\title A note on the mean value of the zeta and 
{\stitle L\/}-functions. XIV} 
\vskip 1cm
\centerline{By Yoichi {\csc Motohashi}${}^{\ast)}$}
\vskip 1cm 
\hsize=16.5cm
{\bf Abstract:} The aim of the present note is to develop a
study on the feasibility of a unified
theory of mean values of automorphic $L$-functions, a desideratum in
the field.  This is an outcome of the investigation commenced with the
part XII ([13]), where a framework was laid on the basis of the theory
of automorphic representations, and a general approach to the mean
values was envisaged. Specifically, it is shown here that the
inner-product method, which was initiated by A. Good [6] and greatly
enhanced by M. Jutila [8], ought to be brought to perfection so that
the mean square of the $L$-function attached to any cusp form on the
upper half-plane could be reached within the notion of automorphy. The
Kirillov map is our key implement. Because of
its geometric nature, our argument appears to extend to bigger linear
Lie groups. This note is essentially self-contained.
\par
\hoffset=-0.2cm
\doublecolumns
\hsize=227.38682pt
\bigskip
\footnote{}{\small 2000 Mathematics Subject 
Classification:  11F70
\par
${}^{\ast)}$ Department of Mathematics, College of Science and
Technology, Nihon University, Surugadai, Tokyo 101-8308.}
{\bf 1. Basic notion.} We collect here basics from the
theory of automorphic representations. In the next section
our problem is made precise. An
idea to deal with it is given in the third section.
\smallskip
Let $G=\r{PSL}(2,{\Bbb
R})$ and $\varGamma=\r{PSL}(2,{\Bbb Z})$. Write
$$
\eqalign{
&\r{n}[x]=\left[\matrix{1&x\cr&1}\right],\quad
\r{a}[y]=\left[\matrix{\sqrt{y}&\cr&1/\sqrt{y}}\right],\cr
&\r{k}[\theta]=\left[\matrix{\phantom{-}\cos\theta&
\sin\theta\cr-\sin\theta &\cos\theta}\right],
}
$$
and 
$$
\eqalign{N&=\left\{\r{n}[x]: x\in\B{R}\right\},\quad
A=\left\{\r{a}[y]: y>0\right\},\cr
K&=\left\{\r{k}[\theta]:\theta\in\B{R}/\pi\B{Z}\right\},
}
$$ 
Thus $G=NAK$ is the Iwasawa decomposition of the Lie group $G$. 
We read it as 
$$
G\ni\r{g}=\r{n}\r{a}\r{k}=\r{n}[x]\r{a}[y]\r{k}[\theta].
$$
The coordinate $(x,y,\theta)$ will retain this definition. 
\par
The center of the universal enveloping algebra 
of $G$ is the polynomial ring on the Casimir element  
$$
\Omega=y^2\left(\partial_x^2+
\partial_y^2\right)-y\partial_x\partial_\theta.
$$
\par
The Haar
measures on the groups
$N$, $A$, $K$, $G$ are normalized, respectively, by $d\r{n}=dx$,
$d\r{a}=dy/y$, $d\r{k}=d\theta/\pi$,
$d\r{g}=d\r{n}d\r{a}d\r{k}/y$,  with Lebesgue measures $dx$,
$dy$, $d\theta$. 
\par
The space $L^2(\varGamma\backslash G)$ is composed of all
left $\varGamma$-automorphic functions on $G$, vectors for short, which
are square integrable over $\varGamma\backslash G$ against $d\r{g}$.
Elements of $G$ act unitarily on vectors from the right. We have the
orthogonal decomposition
$$ 
L^2(\varGamma\backslash G)=\B{C}\cdot1\bigoplus
{}^0\!L^2(\varGamma\backslash G)\bigoplus
{}^e\!L^2(\varGamma\backslash G)
$$ 
into invariant subspaces. Here ${}^0\!L^2$ is the cuspidal subspace
spanned by vectors whose Fourier expansions with respect to the left
action of $N$ have vanishing constant terms. The subspace
${}^e\!L^2$ is spanned by integrals of Eisenstein series. 
Invariant subspaces of $L^2(\varGamma\backslash G)$ and
$\varGamma$-automorphic representations of
$G$ are interchangeable  concepts, and we shall refer to them in a
mixed way. 
\par
The cuspidal subspace decomposes into irredu\-cible subspaces
$$ 
{}^0\!L^2(\varGamma\backslash G)=\overline{\bigoplus V}.
$$ 
The operator $\Omega$ becomes a constant multiplication in each
$V$:
$$
\Omega|_{V^\infty}=\left(\nu^2-{1\over4}\right)\cdot1\,,
$$ 
where $V^\infty$ is the set of all infinitely differentiable 
vectors in $V$. Under our present supposition, $V$ belongs to either
the unitary principal series or the discrete series; accordingly,
we have $\nu\in i\B{R}$ or $\nu=\ell-{1\over2}$ $(1\le\ell\in\B{Z})$. 
\par
The right action of $K$ induces the decomposition of
$V$ into $K$-irreducible subspaces
$$ 
V=\overline{\bigoplus_{p=-\infty}^\infty V_p}\,,\quad \dim V_p\le
1.
$$ 
If it is not trivial, $V_p$ is spanned by a
$\varGamma$-automorphic function on which the right translation by
$\r{k}[\theta]$ becomes the multiplication by the factor
$\exp(2ip\theta)$. It is called a $\varGamma$-automorphic form of
spectral parameter $\nu$ and weight $2p$.
\par
Let $V$ be in the unitary principal series.  Then $\dim V_p=1$ for all
$p\in\B{Z}$ and there exists a complete orthonormal system
$\{\varphi_p\in V_p:\,p\in\B{Z}\}$ of $V$ such that
$$
\varphi_p(\r{g})=\sum_{\scr{n=-\infty}\atop\scr{n\ne0}}
^\infty{\varrho_V(n)\over\sqrt{|n|}}
\e{A}^{\sgn(n)}\phi_p(\r{a}[|n|]\r{g};\nu),
\leqno(1)
$$ 
where $\phi_p(\r{g};\nu)=y^{\nu+{1\over2}}\exp(2ip\theta)$, and
$$
\eqalign{
&\e{A}^\delta\phi_p(\r{g};\nu)\cr
&=\int_{-\infty}^\infty
\exp(-2\pi i\delta x)\phi_p (\r{w}\r{n}[x]\r{g};\nu)dx,
}
$$ 
with $\r{w}=\left[\matrix{&1\cr-1&\cr}\right]$, the Weyl element.
It should be observed that the coefficients
$\varrho_V(n)$ in $(1)$ do not depend on the weight.
In particular, we have the expansion
$$
\leqalignno{
&\varphi_0(\r{g})={2\pi^{{1\over2}+\nu}\over
\Gamma({1\over2}+\nu)}&(2)\cr
&\cdot\sqrt{y}\sum_{\scr{n=-\infty}\atop\scr{n\ne0}}
^\infty\varrho_V(n)\exp(2\pi inx)K_{\nu}(2\pi|n|y),
}
$$ 
where $K_\nu$ is the $K$-Bessel function of
order $\nu$. This is in fact a real analytic cusp form on the
upper half-plane $\{ x+iy: x\in\B{R}, y>0\}$.
\par
Next, let $V$ be in the discrete series, with the
spectral parameter $\ell-{1\over2}$ as above. We have
$$
\hbox{either\quad
$\displaystyle{V=\overline{\bigoplus_{p=\ell}^\infty   V_p}}$
\quad or\quad
$\displaystyle{V=\overline{\bigoplus_{p=-\infty}^{-\ell} V_p}}$}\,,
$$ 
with $\dim V_p=1$, corresponding to the holomorphic and the
antiholomorphic discrete series. The involution
$\r{g}=\r{nak}\mapsto\r{n}^{-1}\r{a}\r{k}^{-1}$ maps one to
the other; thus we may restrict ourselves to the holomorphic case. 
Then we have a complete orthonormal system
$\{\varphi_p\in V_p:\, p\ge\ell\}$ in $V$ such that
$$
\eqalign{
\varphi_p(\r{g})&=\pi^{{1\over2}-\ell}\left({\Gamma(p+\ell)\over
\Gamma(p-\ell+1)}\right)^{1\over2}\cr
&\cdot\sum_{n=1}^\infty
{\varrho_V(n)\over\sqrt{n}}
\e{A}^+\phi_p\left(\r{a}[n]\r{g};\ell-{1\over2}\right),
}
$$ 
In particular, we have
$$
\leqalignno{
\varphi_\ell&(\r{g})=(-1)^\ell{2^{2\ell}\pi^{\ell+{1\over2}}\over
\sqrt{\Gamma(2\ell)}}\exp(2i\ell\theta)&(3)\cr
&\cdot y^\ell
\sum_{n=1}^\infty\varrho_V(n)n^{\ell-{1\over2}}\exp(2\pi in(x+iy)),
}
$$ 
in which the sum is a holomorphic cusp form of
weight $2\ell$ on the upper half-plane. 
\par
It is convenient to put
$$
\hbox{$\psi_V=\varphi_0$\quad or \quad$\varphi_\ell$}\,,
$$
according to the series to which $V$ belongs, with the specification
$(2)$ and $(3)$.
The right Lie derivatives of $\psi_V$ generate the space $V$.
\bigskip
{\bf 2. The problem.} We define the automorphic $L$-function
associated with the irreducible $\varGamma$-automorphic representation
$V$ by
$$
L_V(s)=\sum_{n=1}^\infty \varrho_V(n)n^{-s},
$$
which converges for $\Re s>1$ and continues to an entire function of
polynomial order in any fixed vertical strip. 
\par
The mean square of $L_V$
is the integral
$$
\e{M}_2(L_V;g)=\int_{-\infty}^\infty\left|L_V\left({1\over2}+it\right)
\right|^2g(t)dt\,,
$$
where the weight function
$g$ is assumed, for the sake of simplicity but without much loss
of generality, to be   even, entire, real on $\Bbb R$, and of fast
decay in any fixed horizontal strip. This is an analogue of the
fourth moment of the Riemann zeta-function
$$
\e{M}_4(\zeta;g)=\int_{-\infty}^\infty\left|\zeta
\left({1\over2}+it\right)
\right|^4g(t)dt\,,
$$
since the product of two values of $\zeta$
corresponds to the Eisenstein series in much the same way as $L_V$ does
to $\psi_V$. These quantities have been major subjects in Analytic
Number Theory, as they are indispensable means to reveal the intriguing
nature of the zeta- and $L$-functions.
\par
A. Good [6] was the first to consider $\e{M}_2(L_V;g)$ on the basis of
automorphy. He dealt with the case where $V$ is in the discrete series.
Later the  present author [11] devised an alternative argument
for the same case. He started with the integral
$$
\int_{-\infty}^\infty L_V(u+it)\overline{L_V(\bar{v}+it)}g(t)dt\,,
$$
which is an entire function of $u,\,v$. The non-diagonal part of this
is, in the region of absolute convergence,
$$
\sum_{m=1}^\infty\sum_{n=1}^\infty{
\varrho_V(n)\overline{\varrho_V(n+m)}\over n^u(n+m)^v}\hat{g}\left(
\log\left(1+{m\over n}\right)\right),\leqno(4)
$$
where $\hat{g}$ is the Fourier transform of $g$. Expressing the
$\hat{g}$-factor in terms of a Mellin inversion, one sees that
what is essential is to relate analytically the function
$$
\sum_{n=1}^\infty{
\varrho_V(n)\overline{\varrho_V(n+m)}\over (n+m)^s}\leqno(5)
$$
to the space $V$. For the discrete series this is by no means
difficult. It suffices to consider the inner product
$$
\langle P_m(\cdot;\xi),\,|\psi_V|^2
\rangle_{\varGamma\backslash G/K}\,,\leqno(6)
$$
where $\xi$ is a complex parameter and $P_m$ is the Poin\-car\'e series
of the Selberg type, i.e., the one obtained by replacing the
factor $\tau$ by the constant $1$ that is implicit in $(10)$ below.
Thus the success of the arguments in [6] and [11] is due much to the
fact that the exponential function
$\exp(2\pi i(x+iy))$ appears in the Fourier expansion $(3)$. Having
such a relation between $(5)$ and $(6)$, the spectral decomposition of
$\e{M}_2(L_V;g)$ with $V$ in the discete series reduces to a matter of
technicalities, though the procedure is never straightforward as can
be seen from these two works. 
\par
On the other hand, when $V$ is in the unitary principal series, 
$(6)$ gives only an expression approximating $(5)$ in an involved way.
The difficulty stems from the expansion $(2)$; that is, the presence
of the $K$-Bessel function is the obstacle. Nevertheless, Jutila [8] 
pressed the matter and could develop a deep asymptotic study
of $\e{M}_2(L_V;g)$. A notable merit of his method is in that it is
applicable to the Riemann zeta and any automorphic $L$-functions
equally, even though it does not give complete spectral
decompositions of the mean values.
\par
As to $\e{M}_4(\zeta;g)$, it is in a mixed status. 
An explicit spectral decomposition was established by the present
author [12, Chapter 4]. Recently R.W. Bruggeman and the present
author [3] gave a new proof of it.
Both the arguments are  sharply different from the inner-product
approach mentioned above. In [12, Chapter 4], it is exploited that
the expression corresponding to $(4)$ has the sum of powers of divisors
function $\sigma_\eta$ in place of $\varrho_V$. Ramanujan's expansion
of $\sigma_\eta$ in terms of additive characters and a use of the
functional equation for the Estermann zeta-function transform the
expression into an instance of sums of Kloosterman sums; then
the Kloosterman-spectral sum formula of N.V. Kuznetsov yields the
spectral decomposition. In contrast, the work [3] dispenses with the
use of the sum formula; instead, an approach
via a particular Poincar\'e series on
$\varGamma\backslash G$ is employed. Such a possibility was indicated
already in [10] (see also [12, Section 4.2]), but its realization took
a long time because of the necessity of a drastic change of means. It
was required to employ the representational approach or to move 
from the upper half-plane to the group $G$. More precisely, the
computation of the projection of the Poincar\'e series to an arbitrary
irreducible subspace became the main issue, and it was
accomplished, in [3], only after the authors had been inspired by the
work [5] due to J.W. Cogdell and I. Piatetski-Shapiro.
\par
At any event, both [12, Chapter 4] and [3] depend much on the
arithmetic peculiarity of the function $\sigma_\eta$, and as such it
does not seem to extend to $\e{M}_2(L_V;g)$ with an arbitrary $V$.
\medskip
Thus, there exists a difference among methods 
for mean values of automorphic $L$-functions. It has long been
desired to find a unified way to treat them, indeed since Good's
pioneering works. This is the problem we are dealing with.
\bigskip
{\bf 3. An idea.}
From the above, one might surmise that Jutila's inner-product argument
[8] should be more on the right track than other approaches, because
of its generality. We are going to show that this might be the case.
Namely, albeit a certain restriction is imposed on the
variables $u$, $v$, we shall prove that $(4)$ can be reached via an
extension of
$(6)$, regardless to which series the representation $V$ belongs. Our
idea is to employ the Kirillov map $\e{K}$ to prove a
quasi-analogue of
$(3)$ for any $V$ in the unitary principal series. 
\medskip
To this end we invoke
\medskip
\noindent {\bf Lemma 1.}\quad{\it Let $\nu\in i\B{R}$, and introduce
the Hilbert space
$$ 
U_\nu=\overline{\bigoplus_{p=-\infty}^\infty
\B{C}\phi_p},\quad
\phi_p(\r{g})=\phi_p(\r{g};\nu),
$$ 
equipped with the norm
$$
\Vert\phi\Vert_{U_\nu}=\sqrt{\sum_{p=-\infty}^\infty
|c_p|^2},\quad\phi= \sum_{p=-\infty}^\infty c_p\phi_p.
$$ 
Then 
$$
\e{K}\phi(u)=\e{A}^{\sgn(u)}\phi(\r{a}[|u|])
$$
 is a unitary map from $U_\nu$ onto
$L^2(\B{R}^\times,
\pi^{-1}d^\times)$.
\/}
\medskip
\noindent 
{\it Proof.\/} Here  $\B{R}^\times=\B{R}\backslash\{0\}$
and $d^\times u=du/|u|$, as usual. This seems due originally to A.A.
Kirillov [9] (see also [5, Section 4.2]). A proof of the unitaricity
is given in  [13, Theorem 1], though disguised in the context of
automorphy. The surjectivity is proved in [3, Lemma 4]. As to
analogues for other series of representations, see [3, Section 4].
\hfill $\square$
\medskip
The following assertion is a consequence:
\medskip
\noindent
{\bf Lemma 2.} {\it Let $V$ be in the unitary principal
series. Let $\alpha$ be an arbitrary complex number with
$\Re\alpha$ being positive and sufficiently large.
Then there exists an element $\Phi(\cdot,\alpha)$ in $V$ such that
$$
\leqalignno{
&\Phi(\r{n}[x]\r{a}[y],\alpha)&(7)\cr
&=y^\alpha
\sum_{n=1}^\infty\varrho_V(n)n^{\alpha-{1\over2}}
\exp(2\pi in(x+iy)).
}
$$
}
\par
\noindent
Remark. This brings us to a situation similar to $(3)$. In fact,
we notice a correspondence between $\alpha$ and $\ell$. Our proof
gives a lower bound for $\Re\alpha$, which is, however, not uniform in
$V$.
\smallskip
\noindent
{\it Proof.\/} 
Let $\nu\in i\B{R}$ be the spectral parameter of $V$,
and let $\phi\in U_\nu$ be such that
$$
\e{K}\phi(u)=\cases{u^\alpha\exp(-2\pi u)& for $ u\ge0$,
\cr \hfil 0\hfil & for $u<0$.}\leqno(8)
$$
This is possible, for $\e{K}$ is surjective and the
member on right side is obviously in $L^2({\Bbb R}^\times,
\pi^{-1}d^\times)$. Let
$$
\phi(\r{g})=\sum_{p=-\infty}^\infty a_p\phi_p(\r{g}),\quad
\phi_p(\r{g})=\phi_p(\r{g};\nu),
$$
where $a_p=a_p(\nu,\alpha)$. 
We put
$$
\Phi(\r{g},\alpha)=\sum_{p=-\infty}^\infty a_p\varphi_p(\r{g}),
$$
with $\varphi_p$ as in $(1)$.
We shall later prove briefly that 
$$
\Phi(\r{g},\alpha)=\sum_{\scr{n=-\infty}\atop\scr{n\ne0}}
^\infty{\varrho_V(n)\over\sqrt{|n|}}
\e{A}^{\sgn(n)}\phi(\r{a}[|n|]\r{g}),\leqno(9)
$$
provided $\Re\alpha$ is sufficiently large. This gives $(7)$, since
$$
\e{A}^{\sgn(n)}\phi(\r{a}[|n|]\r{n}[x]\r{a}[y])
=\exp(2\pi inx)\e{K}\phi(ny).
$$
We shall indicate how to prove $(9)$. This is via an
explicit computation of the coefficients $a_p$. The unitaricity of
$\e{K}$ gives
$$
a_p={1\over\pi}\int_0^\infty  u^{\alpha-1}\exp(-2\pi u)
\overline{\e{A}^{+}\phi_p(\r{a}[u])}du\,.
$$
As is well known, the $\e{A}$-factor can be related to the confluent
hypergeometric function (see, e.g., [3, (2.16)]).
Then we use the formula 7.621(3) of [7]. Or one may rather
proceed directly with the present definition. We
find that
$$
\eqalign{
a_p&=(-1)^p 2^{-2\alpha}\pi^{-\nu-\alpha-{1\over2}}\cr
&\cdot{\Gamma\left(\alpha+\nu+{1\over2}\right)
\Gamma\left(\alpha-\nu+{1\over2}\right)\over
\Gamma({1\over2}-\nu+p)\Gamma(\alpha+1-p)}.
}
$$
In particular,
$$
a_p\ll (|p|+1)^{-\Re\alpha-{1\over2}},
$$
as $|p|$ tends to infinity, and $\nu\in i{\Bbb R}$ is bounded.
Thus, indeed $\phi\in U_\nu$ if $\Re\alpha>0$, and $\phi$ becomes
smoother if we take $\Re\alpha$ larger. The confirmation of
$(9)$ follows from this and the uniform bound
$$
\eqalign{
\e{A}^\delta\phi_p(\r{a}[y]&)\ll
(|p|+|\nu|+1)y^{-{1\over2}}\cr
&\cdot\exp\left(-{y\over  
|\nu|+|p|+1}\right)
}
$$ 
(see  [3, (4.5)]). \hfill$\square$
\medskip
Next, we move to an inner-product argument, corresponding to $(6)$: Let
$\tau(\theta)$ be an infinitely differentiable function supported on a
small neighbourhood of $\theta=0$, and
$$
\int_{-{1\over2}\pi}^{{1\over2}\pi}\tau(\theta)d\theta=1.
$$
Let $m$ be a positive integer, and $\Re \xi>1$. Put
$$
h(\r{g})=y^\xi\exp(2\pi mi(x+iy))\tau(\theta).
$$
Further, put
$$
\e{P}h(\r{g})=\sum_{\gamma\in\varGamma_\infty\backslash\varGamma}
h(\gamma\r{g}),\quad \varGamma_\infty=N\cap\varGamma.\leqno(10)
$$
This is in $L^2(\varGamma\backslash G)$. Then, consider the
inner product
$$
\langle \e{P}h,|\Phi|^2\rangle_{\varGamma\backslash G}.
$$
Let us assume that $\alpha$ is positive and
sufficiently large. The unfolding argument gives
$$
\leqalignno{
&\langle\e{P}h,|\Phi|^2\rangle_{\varGamma\backslash G}\cr
&={1\over\pi}
\int_0^\infty\int_0^1y^{\xi-2}\exp(2\pi mi(x+iy))\cr
&\cdot\int_{-{1\over2}\pi}^{{1\over2}\pi}\tau(\theta)|\Phi(\r{n}[x]
\r{a}[y]\r{k}[\theta],\alpha)|^2d\theta dxdy.
}
$$
Thus
$$
\eqalign{
\lim_\tau\, \langle \e{P}h&,|\Phi|^2
\rangle_{\varGamma\backslash G}\cr
={1\over\pi}&
\int_0^\infty\int_0^1y^{\xi-2}\exp(2\pi mi(x+iy))\cr
&\cdot|\Phi(\r{n}[x]\r{a}[y],\alpha)|^2dxdy,
}
$$
where the support of $\tau$ shrinks to the point $0$. The expression
$(7)$ implies that
$$
\eqalign{
\sum_{n=1}^\infty& 
{\varrho_V(n)\overline{\varrho_V(n+m)}\over
(n+m)^\xi(1+m/n)^{\alpha-{1\over2}}}\cr
&={\pi(4\pi)^{\xi+2\alpha-1}
\over \Gamma(\xi+2\alpha-1)}
\lim_\tau\, \langle \e{P}h,|\Phi|^2
\rangle_{\varGamma\backslash G}\,.
}
$$
\par
With this, we may use the
argument of [11, Section 1] and attain the inner sum of $(4)$. In fact,
it suffices for us to multiply both sides by the factor
$$
\eqalign{
&m^{-u-v+\xi}\Gamma(u+v-\xi)\cr
&\cdot{1\over2\pi i}\int_{\Im t=-c}
{\Gamma({3\over2}-u-\alpha+it)
\over\Gamma(v+{3\over2}-\alpha-\xi+it)}g(t)dt
}
$$
with $c>0$ sufficiently large,
and integrate with respect to $\xi$ along an appropriate vertical
line. Provided $\alpha$ is sufficiently large and $\Re(u+v)>\Re\xi>1$, 
the necessary absolute convergence holds throughout our procedure. 
Inserting the thus obtained expression into $(4)$, we find that
$(4)$ admits an expression in terms of $\langle \e{P}h,|\Phi|^2
\rangle_{\varGamma\backslash G}$, provided $\Re(u+v)>2$.
\par
This ends the treatment of the unitary principal series. The case of
the Eisenstein series or that pertaining to $\e{M}_4(\zeta;g)$ is
obviously analogous.
\medskip
Therefore we have proved that $(4)$ with $V$ in any series of
representations can be reached within the notion of automorphy,
provided $u$, $v$ are to be restricted appropriately. Admittedly,
it remains for us to discuss the spectral decomposition that should
follow, especially its analytic continuation to the central point
$(u,v) =\left({1\over2},{1\over2}\right)$. Nevertheless, we may
claim that the above supports the view that there ought to exist a
unified theory of mean values of automorphic $L$-functions. 
\medskip
{\bf Concluding remark.} Probably our choice $(7)$, i.e., $(8)$, will
turn out too drastic in practice. We shall then need to take into
account the smoothness of $\Phi$ when the variable $\r{g}$ approaches
to the Bruhat cell $NA$ from inside the big cell. Namely, we expect
that we shall have to use instead a sequence whose limit is the
present $\phi$, in a way similar to the situation that is experienced
in [3] with the seed function of the Poincar\'e series used there.
In this context, the above should be regarded as a precursor of a
more rigorous discussion to come.
\par
Our argument can readily be extended to
the setting
$G=\r{PSL}(2,\B{C})$ and
$\varGamma=\r{PSL}(2,\B{Z}[\sqrt{-1}])$. All necessary facts are given
in [1, Part XIII] and [2] (see also [1, Part X]). Bigger groups could
also be taken into consideration. For instance, we
presume that a certain {\it double} mean value of the 12th
power of the Riemann zeta-function could be grasped within the setting
$G=\r{PSL}(3,\B{R})$ and $\varGamma=\r{PSL}(3,\B{Z})$. Here the term
double transpires from the real rank of $G$. D. Bump's work [4]
is relevant to this motivation of ours.
\medskip
{\csc Acknowledgement.} We are much indebted to R.W. Bruggeman for his
expert comments to an earlier version of the present note. We 
thank to M. Jutila and A. Ivi\'c for their encouraging comments.

\bigskip
\centerline{\bf References}
\medskip
\item{[1]}  Bruggeman, R.W., and Motohashi, Y.:  A note
on the mean value of the zeta and $L$-functions, X. Proc.\ Japan
Acad., {\bf 77A}, 111--114 (2001); XIII, ibid {\bf 78A},
87--91 (2002).
\item{[2]}  Bruggeman, R.W., and Motohashi, Y.: Sum
formula   for Kloosterman sums and the fourth moment of the Dedekind
zeta-function   over the Gaussian number field. Functiones et
Approximatio, {\bf 31}, 7--76 (2003).
\item{[3]} Bruggeman, R.W., and Motohashi, Y.: A new approach to the
spectral theory of the fourth moment of the Riemann zeta-function.
Submitted.
\item{[4]} Bump, D.: Automorphic Forms on $\r{SL}(3,\B{R})$.
Springer Lect.\ Notes in Math., {\bf 1083}, Springer-Verlag,
Berlin etc., pp.\ 1--184 (1984).
\item{[5]} Cogdell, J.W., and Piatetski-Shapiro, I.:  The
Arithmetic and Spectral Analysis of Poin\-car\'e Series. Academic
Press, San Diego, pp.\ 1-- 192 (1990).
\item{[6]} Good, A.: Beitraege zur Theorie der Dirichletreihen die
Spitzenformen zugeordenet sind. J. Number Theory, {\bf 13}, 18--65
(1981).
\item{[7]} Gradshteyn, I.S., and Ryzhik, I.M.: Tables of
Integrals, Series and Products.  Academic Press, San Diego, 
pp.\ 1--1160 (1979).
\item{[8]} Jutila, M.: Mean values of Dirichlet series via  
Laplace transforms. In: Analytic Number Theory, Proc.\ 39th
Taniguchi Intern.\ Symp.\ Math., Kyoto 1996, ed.\ Y. Motohashi,
Cambridge Univ.\ Press, Cambridge, pp.\ 169--207 (1997).
\item{[9]} Kirillov, A.A.: On $\infty$-dimensional unitary
representations of the group of second-order matrices with elements
from a locally compact field. Soviet Math.\ Dokl., {\bf 4},
748--752 (1963).
\item{[10]} Motohashi, Y.: The fourth power mean of the
Riemann zeta-function. In: Proc.\ Conf.\ Analytic Number Theory,
Amalfi 1989, eds.\ E. Bom\-bieri et al., Univ.\ di Salerno, Salerno,
pp.\ 325--344 (1992).
\item{[11]}  Motohashi, Y.: The mean square of Hecke
$L$-series attached to holomorphic cusp-forms. Ko\-kyuroku Res.\
Inst.\ Math., Kyoto Univ., {\bf 886}, 214--227 (1994).
\item{[12]} Motohashi, Y.:  Spectral Theory of the Riemann
Zeta-Function. Cambridge Univ.\ Press, Cambridge, pp.\ 1--228
(1997).
\item{[13]} Motohashi,Y.: A note on the mean value of the zeta and
$L$-functions, XII. Proc.\ Japan Acad., {\bf 78A}, 36--41 (2002).

\bye